\documentclass{amsart}

\usepackage{amscd}
\usepackage{amssymb}

\newtheorem{theo}{Theorem}[section]
\newtheorem{prop}[theo]{Proposition}

\newtheorem{lemm}[theo]{Lemma}
\newtheorem{defi}[theo]{Definition}

\title{QUANTUM AUTOMORPHISM GROUPS OF FINITE GRAPHS}
\date{} 
\author{Julien Bichon}

\makeatletter
\renewcommand{\@makefnmark}{}
\makeatother

\begin{document}

\begin{abstract}
We determine the quantum automorphism groups of finite graphs.
These are quantum subgroups of the quantum permutation groups defined by
Wang. The quantum automorphism group is a stronger invariant for
finite graphs than the usual automorphism group.
We get a quantum dihedral group $D_4$.
\end{abstract}

\maketitle


\section{Introduction}
The problem of defining the quantum automorphism group of a quantum
space arises quite naturally in noncommutative geometry (\cite{[Co],[Ma]}).
In his recent paper \cite{[Wa]}, S. Wang solves this problem
for the finite space $X_n$ with $n$ points. The main step
was to formulate the universal problem the quantum automorphism group
must solve. The constructed object $A_{aut}(X_n)$ is a compact quantum
group in the sense of Woronowicz (\cite{[W1]}) and is called ``the
quantum permutation group on $n$ symbols''.
Loosely speaking $A_{aut}(X_n)$ is the $C^*$-algebra of functions
on the usual permutation group where the commutativity relations
have been forgotten (and is infinite-dimensional for $n \geq 4$).
Wang also discusses the quantum automorphism groups of 
noncommutative finite-dimensional $C^*$-algebras (see the remark at
the end of section 2). The representation theory of $A_{aut}(X_n)$
has been described by \ T. Banica very recently (\cite{[Ba2]}): the irreducible
representations of $A_{aut}(X_n)$ have the same fusion rules
as the ones of $SO(3)$ (if $n \geq 4$).

\medskip

In this paper we find the quantum automorphism groups of finite graphs.
These are quantum subgroups of the quantum permutation groups.
Our results clearly illustrate the richness of Wang's principle
for quantum automorphism groups.
For example let us consider the graph with $n$ vertices and without edges:
the quantum automorphism group is obviously the quantum permutation 
group. Now let us consider the complete graph with $n$ vertices
where every pair of vertices is an edge: the quantum automorphism group
is the {\sl usual} permutation group. This example shows that
that the construction of the quantum automorphism group is far more
involved than a simple freeness procedure from the usual
automorphism group to the quantum one. It also shows
that the quantum automorphism group is a stronger invariant for
finite graphs than the usual one. This leaves an open door for applications
of the same ideas
to more sophisticated structures in category theory or algebraic topology.

For a special graph with 4 vertices and 4 edges (which is not a square)
we obtain a quantum dihedral group $D_4$ whose $C^*$-algebra
is infinite-dimensional.

\medskip

Our work is organized as follows. In section 2 we recall
some basic definitions and results on compact quantum groups
and quantum automorphism groups and in section 3 we describe the 
quantum automorphism groups of finite graphs.

\medskip

\noindent
{\bf Notations}. All algebras will be unital.

Let  $A = (A,m,u,\Delta,\varepsilon,S)$ be a (complex) Hopf algebra.
 The multiplication will be denoted by  $m$,   $u : \mathbb C
\rightarrow A$ is the unit of $A$, while $\Delta$, $\varepsilon$ and $S$
are respectively the  comultiplication, the counit and the antipode of $A$.

When $A$ is a $*$-algebra and $u \in M_n(A)$ is a matrix, the matrix
$\overline u$ is the matrix $(\overline u)_{ij} = u_{ij}^*$ while 
the transpose matrix of $u$ is denoted by $^tu$
and $u^* = {^t \overline u}$.

\section{Compact quantum groups and quantum automorphism groups}

We first recall some basic notions and results which will
be used freely in the rest of the paper.

\smallskip

Let us recall that a {\bf compact quantum group} (\cite{[W1],[W3]}) is a
pair $(A, \Delta)$ where $A$ is a $C^*$-algebra (with unit)
and $\Delta : A \longrightarrow A \otimes A$ is a coassociative
$*$-homomorphism such that the sets
$\Delta(A)(A\otimes 1)$ and $\Delta(A)(1 \otimes A)$ are both
dense in $A\otimes A$. By abuse of notation a compact quantum group
is often identified with its underlying $C^*$-algebra.
A {\bf morphism} between compact quantum groups $A$ and $B$
is a $*$-homomorphism $\pi : A \longrightarrow B$
such that $\Delta \circ \pi = (\pi \otimes \pi) \circ \Delta$ (\cite{[W]}).

\smallskip

Given a compact quantum group $A$, there is a canonically defined Hopf
$*$-algebra $A^o$ (which we call {\bf the algebra of representative
functions}) which is dense in $A$. The existence of the Haar measure
on $A$ (\cite{[V],[W3]}) shows that furthermore $A^o$ is {\bf unitarizable}
(\cite{[W1]}): for every matrix $u = (u_{ij}) \in M_n(A)$ such that
$\Delta(u_{ij}) = \sum u_{ik} \otimes u_{kj}$ and $\varepsilon(u_{ij})
= \delta_{ij}$, there exists a matrix $F \in GL_n(\mathbb C)$ such that
the matrix $FuF^{-1}$ is unitary (in other words every $A^o$-comodule is
unitarizable).
If $\pi : A \longrightarrow B$ is a morphism of compact quantum groups then
$\pi(A^o) \subset B^o$ and $\pi$ is a Hopf algebra morphism (\cite{[W]}).

\smallskip

Conversely if $A^o$ is a unitarizable Hopf $*$-algebra (a CQG algebra
in \cite{[DK]}), the upperbound of $C^*$-semi-norms exists on
$A^o$ (since $A^o$ is generated by the entries of unitary matrices)
and is a $C^*$-norm (use the regular representation, see \cite{[DK]}, 4.4).
Let $C^*(A^o)$ be the enveloping $C^*$-algebra of $A^o$.
Then $C^*(A^o)$, endowed with the obvious coproduct, is a compact quantum
group. This construction is often useful when one deals with
universal problems (\cite{[W2]} or \cite{[Wa]}). We will use 
it since it is more precise than the direct
$C^*$-construction for the definitions by generators and relations
(no additional relations from $C^*$-norms).   

When $A^o$ is a matrix Hopf $*$-algebra, ie is generated as a $*$-algebra
by entries $u_{ij}$ of a matrix  
$u = (u_{ij}) \in M_n(A^o)$ such that
$\Delta(u_{ij}) = \sum u_{ik} \otimes u_{kj}$ and $\varepsilon(u_{ij})
= \delta_{ij}$ (this condition is equivalent for $A^o$ to be a
finite-type $*$-algebra), there is an easy way to see if  
$A^o$ is unitarizable: indeed $A^o$ is unitarizable if and only if
there are matrices $F \in GL(\mathbb C)$  and
$G \in GL_n(\mathbb C)$ such that
the matrices $FuF^{-1}$ and $G \overline u G^{-1}$ are unitary,
see \cite{[DK]}, 2.4.

\medskip

An {\bf action} of a compact quantum group $A$ on a $C^*$-algebra $Z$ is 
a unital $*$-homomorphism $\alpha : Z \longrightarrow Z \otimes A$
such that there is a dense sub-$*$-algebra $Z^o$ of $Z$
for which $\alpha$ restricts to a coaction 
$\alpha : Z^o \longrightarrow Z^o \otimes A^o$ : ie $Z^o$ is a right
$A^o$-comodule algebra.
The {\bf category of compact quantum transformation groups of $Z$} 
is the category 
whose objects are compact quantum groups acting on $Z$ and whose
morphisms are coaction preserving morphisms of compact quantum
groups (see \cite{[Wa]} for more details).

Let $X$ be a space. The automorphism group of $X$ is a final
object is the category of transformation groups of $X$. This observation
leads to the following definition:

\begin{defi}
Let $Z$ be a $C^*$-algebra. The {\bf quantum automorphism
group of  $Z$} is a compact quantum group $A$ acting on 
$Z$ by $\alpha : Z \longrightarrow Z \otimes A$ and satisfying the 
following universal property:

If $B$ is a compact quantum group acting on $Z$ by 
$\beta : Z \longrightarrow Z \otimes A$, there is a unique 
compact quantum group morphism $\pi : A \longrightarrow B$ such
that $(1 \otimes \pi) \circ \alpha = \beta$ 

When $X$ is a compact space, the quantum automorphism
group of $X$ is the quantum automorphism group
of the $C^*$-algebra $C(X)$. 
\end{defi}

\noindent
{\bf Remarks}.1) When $Z$ is a finite-dimensional $C^*$-algebra, it is 
natural to think that the quantum automorphism group
must be compact. For general $C^*$-algebras,
the object defined above should be called the compact
quantum automorphism group.

\noindent
2) It is clear that the quantum automorphism group, if it exists,
is unique up to isomorphism. It may not exist (see \cite{[Wa]}).

\medskip

Let  $n \in {\mathbb N^*}$ and let $X_n$ be the space with $n$ points.
The $C^*$-algebra $C(X_n)$ is the universal $*$-algebra
with generators $(e_i)_{1 \leq i \leq n}$ and relations
$e_i^* = e_i$ ; $e_i e_j = \delta_{ij} e_i$  and 
$\sum_{i=1}^ne_i = 1$ ; $1 \leq i,j \leq n$. 
The following theorem gives the first example of
quantum automorphism group:

\begin{theo} (\cite{[Wa]}, 3.1)
Let $A_{aut}^o(X_n)$ be the universal (complex) algebra with
generators  $(x_{ij})_{1\leq i,j \leq n}$ and with relations:
$$x_{ij} x_{ik} = \delta_{jk} x_{ij} \ ; \
x_{ji} x_{ki} = \delta_{jk} x_{ji} \ ; \
\sum_l x_{il} = 1 = \sum_l x_{li} \ 
; \ 1 \leq i,j,k \leq n$$
1) There is a Hopf $*$-algebra structure on $A_{aut}^o(X_n)$ defined
by:
$$ x_{ij}^* = x_{ij} \ ; \ 
\Delta(x_{ij}) = \sum_k x_{ik} \otimes x_{kj} \ ; \
\varepsilon(x_{ij}) = \delta_{ij} \ ; \
S(x_{ij}) = x_{ji} \ ;  \ 1\leq i,j \leq n$$
Furthermore $A_{aut}^o(X_n)$ is a unitarizable Hopf $*$-algebra.

\smallskip

\noindent
2) Let $A_{aut}(X_n)$ be the enveloping $C^*$-algebra of
$A_{aut}^o(X_n)$. There is an action of $A_{aut}(X_n)$ on
$C(X_n)$ defined by $\alpha(e_i) = \sum_j e_j \otimes x_{ji}$
and $A_{aut}(X_n)$ is the quantum automorphism
group of the space $X_n$. The spectrum of $A_{aut}(X_n)$ is the 
symmetric group $S_n$.
\end{theo} 

The above quantum group is called {\bf the quantum permutation group} on
$n$ symbols. When $n \geq 4$ then $A_{aut}(X_n)$ is a noncommutative
and infinite-dimensional $C^*$-algebra. 

\medskip

\noindent
{\bf Remark}. Let $Z$ be a finite-dimensional noncommutative $C^*$-algebra.
Wang shows in \cite{[Wa]} that the quantum automorphism of $Z$
does not exist. However if $\psi = {\rm Tr}$ is a trace on $Z$,
he shows the existence of the quantum automorphism group of the pair
$(Z,\psi)$ (see \cite{[Wa]} for the definition). We note here
that this result can be obtained in the algebraic category
of quantum transformation groups of $(Z,\psi)$, and the Hopf algebra
representing the algebraic quantum automorphism group is 
a unitarizable Hopf $*$-algebra whose enveloping $C^*$-algebra
is isomorphic with the quantum group described in theorem
5.1 of \cite{[Wa]}.

\section{Quantum automorphism group of a finite graph}

Let us recall that a {\bf finite graph} $\mathcal G = (V,E)$
consists of a finite set of vertices $V$ and a set of
edges $V \subset E \times E$.

Let $s : E \rightarrow V$ (resp. $t : E \rightarrow V$) be the
source map (resp. the target map). The source and target maps induce
$*$-homomorphisms $s_*, t_* : C(V) \longrightarrow C(E)$.

An (usual) automorphism of a graph is a permutation of the vertices
which preserves the set of edges.

\begin{defi}
An {\bf action of a compact quantum group $A$ on a finite graph
} $\mathcal G = (V,E)$ consists of an action of $A$ on the
set of vertices
$\alpha : C(V) \longrightarrow C(V) \otimes A$ and an action of 
$A$ on the set of edges $\beta : C(E) \longrightarrow C(E) \otimes A$
such that the following diagram $(\star)$ commutes:

\begin{equation*}
\begin{CD}
C(V) \otimes C(V) @> \alpha \underline \otimes \alpha >>
C(V) \otimes C(V) \otimes A \\
@V{m \circ (s_* \otimes t_*)}VV @VV{(m \circ (s_* \otimes t_*))\otimes id_A}V\\
C(E) @> \beta >> C(E) \otimes A
\end{CD}
\end{equation*}
where $m : C(E) \otimes C(E) \rightarrow C(E)$ is the multiplication
map of $C(E)$ and $\alpha \underline \otimes \alpha$ is the tensor
product of the coaction $\alpha$ by itself.

The {\bf quantum automorphism group of a finite graph $\mathcal G = (V,E)$}
is a compact quantum group $A$ acting on $\mathcal G$ by 
$\alpha : C(V) \longrightarrow C(V) \otimes A$ and 
$\beta : C(E) \longrightarrow C(E) \otimes A$, and satisfying
the following universal property:

If $B$ is a compact quantum group acting on $\mathcal G$ by
$\alpha' : C(V) \longrightarrow C(V) \otimes B$ and 
$\beta' : C(E) \longrightarrow C(E) \otimes B$, there is
a unique morphism of compact quantum groups $\phi : A \longrightarrow B$
such that $(id \otimes \phi) \circ \alpha = \alpha'$ and
 $(id \otimes \phi) \circ \beta = \beta'$.
\end{defi}  

It is easy to see that this definition coincide
with the usual one for compact groups. In the group case the action
on the edges is entirely determined by the action on
the vertices. This is also true in the quantum case (since 
the map $m \circ (s_* \otimes t_*)$ is surjective).

\medskip

\noindent
{\bf Notations}. Let $\mathcal G = (V,E)$ be a finite graph with
$n$ vertices $V = \{1,...,n\}$ and $m$ edges $E = \{\gamma_1,..., \gamma_m\}$.
Let $(e_i)_{1 \leq i \leq n}$ be the elements of $C(V)$
defined by $e_i(k) = \delta_{ik}$, $1 \leq k \leq n$ and let 
$(f_j)_{1\leq j \leq m}$ be the elements of $C(E)$ defined
by $f_j(\gamma_l) = \delta_{jl}$, $1\leq l \leq m$.

\begin{theo}
Let $\mathcal G = (V,E)$ be a finite graph with $n$ vertices and
$m$ edges: $E = \{\gamma_1,..., \gamma_m\}$.
Let $A_{aut}^o(\mathcal G)$ be the universal complex algebra
with generators $(X_{ij})_{1\leq i,j \leq n}$ and relations:
$$ 
X_{ij} X_{ik} = \delta_{jk} X_{ij} \  ; \ 
X_{ji} X_{ki} = \delta_{jk} X_{ji} \  ; \ 
\sum_{l=1}^n X_{il} = 1 = \sum_{l=1}^n X_{li} \  
, \ 1 \leq i,j,k \leq n \leqno(3.1)$$
$$ 
X_{s(\gamma_j)i} X_{t(\gamma_j)k} = X_{t(\gamma_j)k} X_{s(\gamma_j)i} = 0 
\leqno(3.2)$$
$$
\quad \quad \quad \quad \quad \quad
X_{i s(\gamma_j)} X_{k t(\gamma_j)} = X_{k t(\gamma_j)} X_{i s(\gamma_j)}
=0 \quad, \quad \gamma_j \in E, \ (i,k) \not\in E$$
$$ 
X_{s(\gamma_j)s(\gamma_l)} X_{t(\gamma_j)t(\gamma_l)} =
X_{t(\gamma_j)t(\gamma_j)} X_{s(\gamma_j)s(\gamma_l)}
\quad , \ \gamma_j,\gamma_l \in E \leqno(3.3)$$
$$ 
\sum_{l=1}^m X_{s(\gamma_l)s(\gamma_j)} X_{t(\gamma_l)t(\gamma_j)}
= 1 = \sum_{l=1}^m X_{s(\gamma_j)s(\gamma_l)} X_{t(\gamma_j)t(\gamma_l)}
\quad , \quad \gamma_j \in E \leqno(3.4) $$
1) There is a Hopf $*$-algebra structure on  $A_{aut}^o(\mathcal G)$
defined by:
$$ X_{ij}^* = X_{ij} \ ; \ 
\Delta(X_{ij}) = \sum_{k=1}^n X_{ik} \otimes X_{kj} \ ; \
\varepsilon(X_{ij}) = \delta_{ij} \ ; \
S(X_{ij}) = X_{ji} \ ; \ 1\leq i,j \leq n$$
and $A_{aut}^o(\mathcal G)$ is a unitarizable Hopf $*$-algebra.

\noindent
2) Let $A_{aut}(\mathcal G)$ be the enveloping $C^*$-algebra 
of $A_{aut}^o(\mathcal G)$. The formulas:
$$ \alpha(e_i) = \sum_{k=1}^n e_k \otimes X_{ki} \quad, \quad 1 \leq i \leq n$$
$$ \beta(f_j) = \sum_{l=1}^m f_l \otimes X_{s(\gamma_l)s(\gamma_j)}
X_{t(\gamma_l)t(\gamma_j)} \quad , \quad 1\leq j \leq m $$
define an action of $A_{aut}(\mathcal G)$ on $\mathcal G$
and $A_{aut}(\mathcal G)$ is the quantum automorphism group
of $\mathcal G$. The spectrum of $A_{aut}(\mathcal G)$ is the
usual automorphism group of $\mathcal G$.
\end{theo}

\noindent
{\bf Proof}. Let $I$ be the two-sided ideal of $A_{aut}^o(X_n)$
generated by relations (3.2), (3.3) and (3.4).
Obviously we have $A_{aut}^o(\mathcal G) \cong A_{aut}^o(X_n)/I$.
It is easy to see that $I$ is a $*$-ideal and hence
$A_{aut}^o(\mathcal G)$ is a $*$-algebra. We must prove that
the maps $\varepsilon$, $\Delta$ and $S$ of the theorem
are well defined. It is easy to check that the character $\varepsilon$
and the anti-homomorphism $S$ are well defined.
Let $\gamma_j \in E$ and $(i,k) \in V \times V$ such that 
$(i,k) \not\in E$. We have
\begin{eqnarray*}
\Delta(X_{s(\gamma_j)i} X_{t(\gamma_j)k}) & = & 
\sum_{p,q}^n X_{s(\gamma_j)p} X_{t(\gamma_j)q} \otimes X_{pi} X_{qk} \\
& = & \sum_{p,q \ (p,q) \in E}
 X_{s(\gamma_j)p} X_{t(\gamma_j)q} \otimes X_{pi} X_{qk} 
\quad ({\rm by \ (3.2)}) \\
& = & \sum_{l=1}^m 
 X_{s(\gamma_j)s(\gamma_l)} X_{t(\gamma_j)t(\gamma_l)} \otimes 
X_{s(\gamma_l)i} X_{t(\gamma_l)k} = 0 \quad ({\rm by \ (3.2)}).
\end{eqnarray*}
In the same way we have
$$\Delta(X_{t(\gamma_j)k} X_{s(\gamma_j)i}) = 
0 = \Delta(X_{is(\gamma_j)}X_{kt(\gamma_j)}) =
\Delta(X_{kt(\gamma_j)} X_{is(\gamma_j)}).$$
Let $\gamma_j$, $\gamma_l \in E$. We have
\begin{eqnarray*}
\Delta(X_{s(\gamma_j)s(\gamma_l)} X_{t(\gamma_j)t(\gamma_l)}) & = &
\sum_{i,k}^n X_{s(\gamma_j)i} X_{t(\gamma_j)k} \otimes
X_{is(\gamma_l)} X_{kt(\gamma_l)} \\
({\rm by \ (3.2)}) \quad & = & \sum_{l' = 1}^m X_{s(\gamma_j)s(\gamma_{l'})} 
X_{t(\gamma_j)t(\gamma_{l'})} \otimes X_{s(\gamma_{l'})s(\gamma_l)}
X_{t(\gamma_{l'})t(\gamma_l)} \quad \\
({\rm by \ (3.3)}) \quad & = & \sum_{l'=1}^{m} X_{t(\gamma_j)t(\gamma_{l'})}
X_{s(\gamma_j)s(\gamma_{l'})} \otimes X_{t(\gamma_{l'})t(\gamma_l)}
X_{s(\gamma_{l'})s(\gamma_l)} \\
({\rm by \ (3.2)}) \quad & = & \sum_{i,k}^n X_{t(\gamma_j)i} X_{s(\gamma_j)k}
\otimes X_{it(\gamma_l)} X_{ks(\gamma_l)} \\
& = & \Delta(X_{t(\gamma_j)t(\gamma_l)} X_{s(\gamma_j)s(\gamma_l)}).
\end{eqnarray*}
Now let $\gamma_j \in V$. We have:
\begin{eqnarray*} 
\Delta(\sum_{l=1}^m X_{s(\gamma_l)s(\gamma_j)} X_{t(\gamma_l)t(\gamma_j)})
& =  & \sum_{l=1}^m \sum_{k,i}^n
X_{s(\gamma_l)i} X_{t(\gamma_l)k} \otimes X_{is(\gamma_j)}
X_{kt(\gamma_j)} \\
({\rm by \ (3.2)}) & = & \sum_{l,l'}^m 
X_{s(\gamma_l)s(\gamma_{l'})} X_{t(\gamma_l)t(\gamma_{l'})} \otimes
X_{s(\gamma_{l'})s(\gamma_j)} X_{t(\gamma_{l'})t(\gamma_j)} \\  
& = & 1 \otimes 1 = \Delta(1) \quad ({\rm by \ (3.4)}).
\end{eqnarray*}
In the same way $\Delta(\sum_{l=1}^m 
X_{s(\gamma_j)s(\gamma_l)} X_{t(\gamma_j)t(\gamma_l)}) = \Delta(1)$.
Therefore $\Delta$ is a well defined algebra morphism and in this way
$A_{aut}^o(\mathcal G)$ is a bialgebra. It is clear that
the antihomomorphism $S$ is an antipode for this bialgebra and
hence $A_{aut}^o(\mathcal G)$ is a Hopf $*$-algebra
which is clearly unitarizable.
There is an obvious surjective Hopf $*$-algebra morphism 
$\pi : A_{aut}^o(X_n) \longrightarrow A_{aut}^o(\mathcal G)$ defined
by $\pi(x_{ij}) = X_{ij}$.

\smallskip

\noindent
2) By relations (3.1) and theorem 2.2 we have an action
$\alpha : C(V) \longrightarrow C(V) \otimes A_{aut}(\mathcal G)$ as defined in
the theorem. Let $\gamma_j$, $\gamma_l$ and $\gamma_{j'} \in E$.
By (3.3) we have
$$ (X_{s(\gamma_l)s(\gamma_j)} X_{t(\gamma_l)t(\gamma_j)})^*
= X_{s(\gamma_l)s(\gamma_j)} X_{t(\gamma_l)t(\gamma_j)}$$
and by (3.1), (3.2) and (3.3) we have
$$X_{s(\gamma_l)s(\gamma_j)} X_{t(\gamma_l)t(\gamma_j)}
X_{s(\gamma_l)s(\gamma_{j'})} X_{t(\gamma_l)t(\gamma_{j'})} =
\delta_{jj'} X_{s(\gamma_l)s(\gamma_j)} X_{t(\gamma_l)t(\gamma_j)}.$$
Finally using (3.4) we see that there is a
$*$-homomorphism $\beta : C(E) \longrightarrow C(E) \otimes 
A_{aut}(\mathcal G)$ as defined in the theorem
and $\beta$ is coassociative (see the calculations
in i)). Let us check that the diagram ($\star$) commutes.
Let $i,k \in V$. We have $\alpha \underline \otimes \alpha(e_i \otimes e_k)
= \sum_{p,q}^n e_p \otimes e_q \otimes X_{pi} X_{qk}$.
We also have $s_*(e_i) = \sum_{l,s(\gamma_l)=i} f_l$ and 
$t_*(e_k) = \sum_{l,t(\gamma_l) = k} f_l$. Hence
\begin{eqnarray*}
((m \circ (s_* \otimes t_*)) \otimes id) \circ 
(\alpha \underline \otimes \alpha) & = & 
\sum_{p,q}^n \sum_{l,(p,q) = \gamma_l} 
f_l \otimes X_{pi} X_{qk} \\
& = & \sum_{j=1}^m f_j \otimes X_{s(\gamma_j)i} X_{t(\gamma_j)k}.
\end{eqnarray*}
On the other hand
\begin{eqnarray*}
\beta \circ m \circ (s_* \otimes t_*) (e_i \otimes e_k)
& = & 
\sum_{l,\gamma_l=(i,k)} \ \sum_{j=1}^m
f_j \otimes X_{s(\gamma_j)s(\gamma_l)} X_{t(\gamma_j)t(\gamma_{l})} \\
& = &
\sum_{j=1}^m f_j \otimes X_{s(\gamma_j)i} X_{t(\gamma_j)k}.
\end{eqnarray*}
In this way we have defined an action of $A_{aut}(\mathcal G)$
on the graph $\mathcal G$.

Let $B$ be a compact quantum group acting on $\mathcal G$ by
$\alpha' : C(V) \longrightarrow C(V) \otimes B$
and $\beta' : C(E) \longrightarrow C(E) \otimes B$.
There are elements $(a_{ij})_{1\leq i,j \leq n}$ of $B$
such that $\alpha'(e_i) = \sum_j e_j \otimes a_{ji}$ and by theorem
2.2 these elements satisfy relations (3.1) (the elements $a_{ij}$
belong to the algebra of representative functions $B^o$).
There are also elements $(b_{jl})_{1\leq j,l \leq m}$ of $B^o$
such that $\beta'(f_j) = \sum_l f_l \otimes b_{lj}$ and theorem
2.2 ensures that they satisfy the same relations as (3.1).
By the commutativity of the diagram ($\star$), we have
$b_{lj} = a_{s(\gamma_l)s(\gamma_j)} a_{s(\gamma_l)s(\gamma_j)}$
 and hence there is
a surjective morphism of compact quantum groups 
$\phi_0 : A_{aut}(X_n) \longrightarrow B$ defined by
$\phi_0(x_{ik}) = a_{ik}$. It remains to check that
the $a_{ik}$'s satisfy relations (3.2), (3.3) and (3.4).
We have $b_{lj}^* = b_{lj}$ and hence
$$a_{t(\gamma_l)t(\gamma_j)} a_{s(\gamma_l)s(\gamma_j)} = 
a_{s(\gamma_l)s(\gamma_j)} a_{t(\gamma_l)t(\gamma_j)}
\quad , \quad \gamma_j,\gamma_l \in E.$$
and relations (3.3) hold. We have $\sum_{l=1}^m b_{lj} = 1
= \sum_{l=1}^m b_{jl}, \forall \gamma_j \in E$ and hence 
relations (3.4) are satisfied.
Using once again the commutativity of the diagram ($\star$),
we get $a_{s(\gamma_j)i} a_{t(\gamma_j)k} = 0$ whenever
 $(i,k) \not\in E$. Using the involution and the antipode of $B^o$,
it is easy to see that the other relations in (3.2) are fulfilled.
Thus we have a morphism of quantum groups 
$\phi : A_{aut}(\mathcal G) \longrightarrow B$
such that $(id \otimes \phi) \circ \alpha = \alpha'$
and $(id \otimes \phi) \circ \beta = \beta'$: this morphism
is obviously unique and
$A_{aut}(\mathcal G)$ is the quantum automorphism group of $\mathcal G$.
The last statement follows immediately from the universal property of 
$A_{aut}(\mathcal G)$. $\square$ 

\bigskip

Let us take a look at some examples. First let us consider the graph
$\mathcal G = (V,E)$ with $n$ vertices and $E = \emptyset$.
It is obvious that $A_{aut}(\mathcal G)$ is the quantum permutation
group $A_{aut}(X_n)$. Less trivially, let us consider the
graph $\mathcal G = (V,E)$ with $n$ vertices 
$V=\{1,...,n\}$ and $E = \{(i,i), i \in E \}$. Then it is easy to
check that relations (3.2), (3.3) and (3.4) all follow
from (3.1) and hence $A_{aut}(\mathcal G) = A_{aut}(X_n)$.

Now let us consider the complete graph $\mathcal G = (V,E)$
with $n$ vertices and $E = V \times V$.
Then relations (3.3) imply that $A_{aut}(\mathcal G)$ is a commutative
$C^*$-algebra. Since the spectrum of $A_{aut}(\mathcal G)$ is the
usual automorphism group of $\mathcal G$ (in this case the symmetric
group $S_n$), we have $A_{aut}(\mathcal G) \cong C(S_n)$.
Hence the quantum automorphism group and the usual one coincide in
this example.

These simple examples show that the quantum automorphism
group is a stronger invariant for finite graphs than
the usual automorphism group.

Another elementary example is the following one. Let us consider the
polygonal graph $\mathcal G$ with $n$ vertices 
$V=\{1,...,n\}$ and $E = \{\gamma_1,...,\gamma_n\}$ where
$\gamma_i = (i,i+1)$, $1\leq i \leq n-1$ and $\gamma_n = (n,1)$.
It is easily seen that $A_{aut}(\mathcal G) \cong 
C(\mathbb Z/n \mathbb Z)$, the algebra of functions on the cyclic
group $\mathbb Z/n \mathbb Z$.

We now examine a less trivial example:

\begin{prop}
Let $\mathcal G = (V,E)$ be the graph with 4 vertices 
$V = \{1,2,3,4 \}$ and 4 edges $E=\{\gamma_1,\gamma_2,\gamma_3,\gamma_4 \}$
where $\gamma_1 = (1,2)$ , $\gamma_2 = (2,1)$, $\gamma_3 = (3,4)$ and
$\gamma_4 = (4,3)$. Then $A_{aut}(\mathcal G)$ is a noncommutative
infinite-dimensional $C^*$-algebra whose spectrum is the dihedral
group $D_4$.
\end{prop}

The quantum group above will be called the {\bf quantum dihedral group} $D_4$.

\medskip

\noindent
{\bf Proof}. It is easily seen that the usual automorphism group
of the above graph is a finite group with 8 elements which is 
the dihedral group $D_4$ and hence by theorem 3.2 the spectrum
of $A_{aut}(\mathcal G)$ is the group $D_4$.

Let us now describe the $*$-algebra $A_{aut}^o(\mathcal G)$.
First let us translate relations (3.2). We get the following relations
(3.2)':

$X_{11} X_{23} = X_{23} X_{11} = 0 = X_{32} X_{11} = X_{11} X_{32},$

$X_{11} X_{24} = X_{24} X_{11} = 0 = X_{42} X_{11} = X_{11} X_{42},$

$X_{12} X_{23} = X_{23} X_{12} = 0 = X_{32} X_{21} = X_{21} X_{32},$

$X_{12} X_{24} = X_{24} X_{12} = 0 = X_{42} X_{21} = X_{21} X_{42},$

$X_{13} X_{21} = X_{21} X_{13} = 0 = X_{31} X_{12} = X_{12} X_{31},$

$X_{13} X_{22} = X_{22} X_{13} = 0 = X_{31} X_{22} = X_{22} X_{31},$

$X_{14} X_{22} = X_{22} X_{14} = 0 = X_{41} X_{22} = X_{22} X_{41},$

$X_{14} X_{21} = X_{21} X_{14} = 0 = X_{12} X_{41} = X_{41} X_{12},$

$X_{31} X_{43} = X_{43} X_{31} = 0 = X_{13} X_{34} = X_{34} X_{13},$

$X_{31} X_{44} = X_{44} X_{31} = 0 = X_{13} X_{44} = X_{44} X_{13},$

$X_{32} X_{43} = X_{43} X_{32} = 0 = X_{23} X_{34} = X_{34} X_{23},$

$X_{32} X_{44} = X_{44} X_{32} = 0 = X_{23} X_{44} = X_{44} X_{23},$

$X_{33} X_{41} = X_{41} X_{33} = 0 = X_{33} X_{14} = X_{14} X_{33},$

$X_{33} X_{42} = X_{42} X_{33} = 0 = X_{33} X_{24} = X_{24} X_{33},$

$X_{34} X_{42} = X_{42} X_{34} = 0 = X_{43} X_{24} = X_{24} X_{43},$

$X_{34} X_{41} = X_{41} X_{34} = 0 = X_{14} X_{43} = X_{43} X_{14}.$

\noindent
Relations (3.3) are translated in the following relations (3.3)':

$X_{11} X_{22} = X_{22} X_{11} \quad , \quad X_{31} X_{42} = X_{42} X_{31},$

$X_{12} X_{21} = X_{21} X_{12} \quad , \quad X_{32} X_{41} = X_{41} X_{32},$

$X_{13} X_{24} = X_{24} X_{13} \quad , \quad X_{33} X_{44} = X_{44} X_{33},$

$X_{14} X_{23} = X_{23} X_{14} \quad , \quad X_{34} X_{43} = X_{43} X_{34}.$

\noindent
Relations (3.4) become relations (3.4)':

\noindent
$X_{11} X_{22} + X_{21} X_{12} + X_{31} X_{42} + X_{41} X_{32} = 1 =
X_{11} X_{22} + X_{12} X_{21} + X_{13} X_{24} + X_{14} X_{23},$

\noindent
$X_{12} X_{21} + X_{22} X_{11} + X_{32} X_{41} + X_{42} X_{31} = 1 =
X_{21} X_{12} + X_{22} X_{11} + X_{23} X_{14} + X_{24} X_{13},$

\noindent
$X_{13} X_{24} + X_{23} X_{14} + X_{33} X_{44} + X_{43} X_{34} = 1 =
X_{31} X_{42} + X_{32} X_{41} + X_{33} X_{44} + X_{34} X_{43},$

\noindent
$X_{14} X_{23} + X_{24} X_{13} + X_{34} X_{43} + X_{44} X_{33} = 1 =
X_{41} X_{32} + X_{42} X_{31} + X_{43} X_{34} + X_{44} X_{33}.$

\noindent
Combining (3.1), (3.2)', (3.3)', (3.4)' we get the following useful
relations:

$X_{11} = X_{22} = X_{11} X_{22} = X_{22} X_{11} \quad , \quad
X_{12} = X_{21} = X_{12} X_{21} = X_{21} X_{12},$

$X_{13} = X_{24} = X_{13} X_{24} = X_{24} X_{13} \quad , \quad
X_{14} = X_{23} = X_{14} X_{23} = X_{23} X_{14},$

$X_{31} = X_{42} = X_{31} X_{42} = X_{42} X_{31} \quad , \quad
X_{32} = X_{41} = X_{32} X_{41} = X_{41} X_{32},$

$X_{33} = X_{44} = X_{33} X_{44} = X_{44} X_{33} \quad , \quad
X_{34} = X_{43} = X_{34} X_{43} = X_{43} X_{34}.$

\smallskip

We are now able to describe the algebra $A_{aut}^o(\mathcal G)$ in a 
simpler way. Let $B^o$ be the universal $*$-algebra with generators
$(y_i)_{1\leq i \leq 8}$ and relations:
$$ y_i^* = y_i^2 = y_i \ , \  1\leq i \leq 8 \quad ;
\quad y_1 y_i = 0 = y_i y_1 \ , \  2\leq i \leq 6  \quad ;$$
$$y_2 y_i = 0 = y_i y_2 \ , \ 3 \leq i \leq 6 \quad ; \quad
y_3 y_i = 0 = y_i y_3 \ , \ i \in \{4,7,8\}  \quad ;$$ 
$$y_4 y_i = 0 = y_i y_4 \ , \ i \in \{7,8\} \quad ; \quad
y_5 y_i = 0 = y_i y_5 \ , \ i \in \{6,7,8\} \quad ;$$
$$y_6 y_i = 0 = y_i y_6 \ , \  i \in \{7,8\} \quad ; \quad
y_7 y_8 = 0 = y_8 y_7.$$
$$y_1 + y_2 + y_3 + y_4 = y_1 + y_2 + y_5 + y_6 = 1 =  y_3 + y_4 + y_7 + y_8
= y_5 + y_6 + y_7 + y_8.$$
The reader will easily check that there is a $*$-homomorphism
$\varphi : B^o \longrightarrow A_{aut}^o(\mathcal G)$ defined by
$\varphi(y_1) = X_{11}$, $\varphi(y_2) = X_{12}$, $\varphi(y_3) = X_{13}$,
$\varphi(y_4) = X_{14}$, $\varphi(y_5) = X_{31}$, $\varphi(y_6) = X_{32}$,
$\varphi(y_7) = X_{33}$, $\varphi(y_8) = X_{34}$.

In the same way there is a $*$-homomorphism
$\psi : A_{aut}^o(\mathcal G) \longrightarrow B^o$ defined by
$\psi(X_{11}) = \psi(X_{22}) = y_1$, $\psi(X_{12}) = \psi(X_{21}) = y_2$,
$\psi(X_{13}) = \psi(X_{24}) = y_3$, $\psi(X_{14}) = \psi(X_{23}) = y_4$,
$\psi(X_{31}) = \psi(X_{42}) = y_5$, $\psi(X_{32}) = \psi(X_{41}) = y_6$,  
$\psi(X_{33}) = \psi(X_{44}) = y_7$, $\psi(X_{34}) = \psi(X_{43}) = y_8$.

We have $\varphi \circ \psi =id$ and $\psi \circ \varphi = id$ and therefore
the $*$-algebras $B^o$ and $A_{aut}^o(\mathcal G)$ are isomorphic.

Let us now show that $B^o$ is noncommutative and infinite-dimensional.
There is a representation $\pi : B^o \longrightarrow M_2(\mathbb C)$
defined by:
$$\pi(y_1) = \left( 
\begin{array}{cc}
1 & 0 \\
0 & 0
\end{array} \right) \ , \ 
\pi(y_2) = \left( 
\begin{array}{cc}
0 & 0 \\
0 & 1
\end{array} \right) , \ 
\pi(y_7) = \left( 
\begin{array}{cc}
1 & 0 \\
1 & 0
\end{array} \right) , \ 
\pi(y_8) = \left( 
\begin{array}{cc}
0 & 0 \\
-1 & 1
\end{array} \right)  
$$
and $\pi(y_3) = \pi(y_4) = \pi(y_5) = \pi(y_6) = 0$. We have
$\pi(y_1y_7) = \pi(y_1)$ while $\pi(y_7 y_1) = \pi(y_7)$: $B^o$ is not commutative.

Let us suppose that $B^o$ is finite-dimensional.
Then $B^o$ would be a finite-dimensional $C^*$-algebra since
$B^o \cong A_{aut}^o(\mathcal G)$. But $\pi(B^o)$ is contained in
the algebra of lower-triangular matrices. This means that the representation
$\pi$ is not semisimple, and hence there is no scalar product on
the vector space $\mathbb C^2$ for which $\pi$ is a $*$-representation.
This a contradiction with the following lemma, which concludes the proof
of proposition 3.3: $A_{aut}(\mathcal G) \cong C^*(B^o)$ is infinite-dimensional. $\square$
  
\begin{lemm}
Let $A$ be a finite-dimensional $C^*$-algebra and let $\pi : A \longrightarrow
{\rm End}(V)$ be a representation of the algebra $A$ on a finite-dimensional
vector space $V$. Then there is a Hilbert space structure on $V$ such
that $\pi$ becomes a $*$-representation.
\end{lemm}

\noindent
{\bf Proof}. Let $G$ be the compact group of unitary elements of $A$.
There is a scalar product on $V$ such that $\pi$ becomes a unitary
representation of $G$. But every element in a $C^*$-algebra is a
linear combination of unitary elements, and hence
$\pi$ is a $*$-representation for the above scalar product. $\square$

\medskip

\noindent
{\bf Remark}. Lemma 3.4 is probably well known. It can be used to give a simple
proof of the existence of the Haar measure on a finite quantum group
in the following way.

Let $A$ be a finite-dimensional Hopf $C^*$-algebra. There is a bijection
(an equivalence of categories) between representations of $A$ and 
$A'$-comodules ($A'$ is the dual Hopf $*$-algebra). An $A'$-comodule $V$
is unitarizable if and only if there is a Hilbert space structure on
$V$ such that the corresponding representation is a $*$-representation.
By lemma 3.4 every $A'$-comodule is unitarizable and hence there is
a positive and faithful Haar measure on $A'$ (see \cite{[DK]}).
The regular representation provides a $C^*$-norm on $A'$ and hence
$A'$ is a Hopf $C^*$-algebra. By duality there is a positive 
and faithful Haar measure on $A$.
In this way we can say that the existence of the Haar measure on a finite
quantum group follows from the existence of the Haar measure on usual compact groups.

\medskip

Proposition 3.3 gives a concrete example of non-trivial quantum subgroup
(ie noncommutative and noncocommutative) of $A_{aut}(X_n)$.
It is expected that the construction described in
\cite{[Ba]} furnishes many other examples of this kind.

\bigskip

\noindent
{\bf Aknowledgements}. The author wishes to thank G. Laffaille
who read the manuscript.

\bigskip

\address{D\'epartement des Sciences Math\'ematiques,
case 051,  Universit\'e Montpellier II.

Place Eug\`ene Bataillon, 34095 Montpellier Cedex 5.}

E-mail : bichon\char64math.univ-montp2.fr

\end{document}